\documentclass [leqno] {amsart}
  
\usepackage {amssymb}
\usepackage {amsthm}
\usepackage {amscd}

\usepackage[all]{xy}
\usepackage {hyperref} 

\begin {document}

\theoremstyle{plain}
\newtheorem{proposition}[subsection]{Proposition}
\newtheorem{lemma}[subsection]{Lemma}
\newtheorem{korollar}[subsection]{Corollary}
\newtheorem{thm}[subsection]{Theorem}
\newtheorem*{thm*}{Theorem}
\newtheorem{conjecture}[subsection]{Conjecture}

\theoremstyle{definition}

\theoremstyle{remark}
\newtheorem{remark}[subsection]{Remark}

\numberwithin{equation}{subsection}

\newcommand{\arir}{\ar@{^{(}->}}
\newcommand{\aril}{\ar@{_{(}->}}
\newcommand{\are}{\ar@{>>}}

\newcommand{\xr}[1] {\xrightarrow{#1}}
\newcommand{\xl}[1] {\xleftarrow{#1}}


\newcommand{\CH}{{\rm CH}}
\newcommand{\Gr}{{\rm Gr}}

\newcommand{\Spec} {{\rm Spec}}

\newcommand{\Hom} {{\rm Hom}}
\newcommand{\End} {{\rm End}}

\newcommand{\Aut}{{\rm Aut}}


\newcommand{\PP} {\mathbb{P}}
\newcommand{\Z} {\mathbb{Z}}
\newcommand{\Q} {\mathbb{Q}}
\newcommand{\C} {\mathbb{C}}

\newcommand{\abs}[1]{\lvert#1\rvert}
\newcommand{\OO}{\mathcal{O}}
\newcommand{\Ql}{\Q_{\ell}}

\title{First coniveau notch of the Dwork family and its mirror}
\author{Andre Chatzistamatiou}
\email{a.chatzistamatiou@uni-due.de}

\begin{abstract}
If $X_{\lambda}$ is a smooth member of the Dwork family over a perfect field $k$, and $Y_{\lambda}$ is its mirror variety, 
then the motives of $X_{\lambda}$ and $Y_{\lambda}$ are equal up to motives that are in coniveau $\geq 1$. 
If $k$ is a finite field, this provides a motivic explanation for Wan's congruence between the 
zeta functions of $X_{\lambda}$ and $Y_{\lambda}$.
\end{abstract}

\maketitle

\section*{Introduction} 

Let $k$ be a field. We consider the Dwork family of hypersurfaces $X_{\lambda}$ in $\PP^n$ defined by 
the equation 
$$
\sum_{i=0}^n X_i^{n+1} + \lambda X_0\dots X_n = 0
$$
with the parameter $\lambda\in k$. The variety $X_{\lambda}$ is a Calabi-Yau manifold when $X_{\lambda}$ is smooth. 
On each member $X_{\lambda}$ there is a group action by the kernel $G$
of the character $\mu_{n+1}^{n+1}\xr{} \mu_{n+1}, (\zeta_i)\mapsto \prod_i \zeta_i,$ given by
$$
G\times X_{\lambda} \xr{} X_{\lambda}, \quad 
(\zeta_0,\dots,\zeta_n)\cdot (x_0:\dots:x_n) = (\zeta_0x_0,\dots,\zeta_nx_n). 
$$
The quotient $X_{\lambda}/G$ is a hypersurface with trivial canonical bundle 
in a toric Fano variety and a singular mirror of $X_{\lambda}$ \cite{B}. 
If $Y_{\lambda}$ is a crepant resolution of $X_{\lambda}/G$ (which exists but in general is not unique) then $(X_{\lambda},Y_{\lambda})$ provides an example of a mirror pair.
Since the birational geometry of $Y_{\lambda}$ is independent of the choice of the resolution a
natural question arises: to compare the birational motives of $X_{\lambda}$ and $Y_{\lambda}$.
For a finite field $k=\mathbb{F}_q$ the number of $\mathbb{F}_{q^m}$-rational points modulo
$q^m$ is a birational invariant and
D.~Wan asked to compare the number of rational points of a mirror pair \cite{W}. 
In the case of the Dwork family he  
proved a mirror congruence formula \cite[Theorem~1.1]{W}:
$$
\# X_{\lambda}(\mathbb{F}_{q^m}) = \# Y_{\lambda}(\mathbb{F}_{q^m}) \quad \text{mod $q^m$}
$$
for every positive integer $m$. Fu and Wan studied more general mirror pairs which come 
from quotient constructions and obtained under certain assumptions on the action of $G$
(see Theorem \ref{thmFW}) a congruence formula \cite{FW}:
\begin{equation}\label{congXXG}
\#X(\mathbb{F}_{q^m})=\#(X/G)(\mathbb{F}_{q^m}) \mod q^m. 
\end{equation}
The same formula is proved in {\cite[Corollary~6.12]{BBE}} with different assumptions.

The purpose of this paper is twofold. The first theorem compares the motives of $X_{\lambda}$
and $Y_{\lambda}$ when $X_{\lambda}$ is a member of the Dwork family, and provides Wan's
congruence formula as a consequence. We also explain what can be expected for general quotient 
constructions in \S 3. In the second theorem we prove a congruence formula for a quotient 
singularity $X/G$ and a resolution of singularities $Y\xr{} X/G$:
$$
\#X/G(\mathbb{F}_{q^m})=\#Y(\mathbb{F}_{q^m}) \quad \mod q^m.
$$ 
Thus \ref{congXXG} is sufficient in order to get 
$\#X(\mathbb{F}_{q^m})=\#Y(\mathbb{F}_{q^m})$ modulo $q^m$.  

We state now our theorems and several consequences. 
By a motive we understand a pair $(X,P)$ with $X$ a smooth projective variety  
and $P\in \CH^{\dim X}(X\times X)\otimes \Q$ a projector. The morphisms are correspondences in rational 
coefficients; the Lefschetz motive is denoted by $\Q(-1):=(\PP^1,\PP^1\times p)$ with $p\in \PP^1(k)$. 
For $X_{\lambda}$  the cycle $P=1/\abs{G} \sum_{g\in G}\Gamma(g),$ where $\Gamma$ denotes the graph, is a projector.  

\begin{thm*}
Let $k$ be a perfect field, and $n\geq 2$. 
We assume that ${\rm char}(k)\nmid n+1$ if the characteristic of $k$ is positive. 
Let $X_{\lambda}$ be a smooth member of the Dwork family. 
Then there are motives $N,N'$ such that
$$
(X_{\lambda},id) \cong (X_{\lambda},P) \oplus N\otimes \Q(-1) \quad \text{and} \quad  (Y_{\lambda},id) \cong (X_{\lambda},P) \oplus N'\otimes \Q(-1).
$$
\end{thm*}

For a finite field $k=\mathbb{F}_q$ the eigenvalues of the geometric Frobenius  acting on 
$H^*_{\text{\'et}}(N\otimes \Q(-1))=H^*_{\text{\'et}}(N)\otimes \Q_l(-1)$ lie in $q\cdot \bar{\Z}$, and by using 
Grothendieck's trace formula this implies Wan's theorem \cite[Theorem~1.1]{W}. For $k=\C$ the theorem of  
Arapura-Kang on the functoriality of the coniveau filtration $N^*$ allows us to conclude that
$$
{\rm gr}^0_{N^*}(H^*(X_{\lambda},\Q)) \cong {\rm gr}^0_{N^*}(H^*(Y_{\lambda},\Q)) 
$$ 
as Hodge structures (see Corollary \ref{consequences}). 
  
We now describe our method. We use birational motives in order to reduce to a statement for 
zero cycles over $\C$: $\CH_0(X_{\lambda})=P\circ \CH_0(X_{\lambda})$, i.e.~$P$ acts as identity. 
To prove this we consider, additionally to $G$, the action of the symmetric group $S_{n+1}$ acting via 
permutation of the homogeneous coordinates. The transpositions act as $-1$ on $H^0(X_{\lambda},\omega_{X_{\lambda}})$
and the quotients $X_{\lambda}/H$ for suitable subgroups $H$ of $G\rtimes S_{n+1}$ can be shown to be $\Q$-Fano 
varieties. By the theorem of Zhang \cite{Z} these are rationally chain connected, which yields sufficiently 
many relations for the zero cycles on $X_{\lambda}$ to prove the claim.      

\begin{thm*}
Let $X$ be a smooth projective $\mathbb{F}_q$-variety with an action of a finite group $G$.
Let $\pi:X\xr{} X/G$ be the quotient, and $f:Y\xr{} X/G$ be a birational map, where $Y$ is 
a smooth projective variety. Then
$$
\#Y(\mathbb{F}_q)=\#X/G(\mathbb{F}_q) \quad \mod q.
$$
\end{thm*}

For the proof we use the action of the geometric Frobenius $F$ on \'etale cohomology. 
Suppose that $Z\subset X/G$ is the set where $f$ is not an isomorphism, then $F$ acts
on the cohomology with support in $Z$ with eigenvalues in $q\bar{\Z}$. This is proved
by reduction to the case $\pi^{-1}(Z)\subset X$ via a trace map argument. Counting 
points with Grothendieck's trace formula yields the result.

\subsection*{Acknowledgments}
I thank Y.~Andr\'e for drawing my attention to D.~Wan's work and for helpful suggestions. 
I thank D.C.~Cisinski, F.~D\'eglise and K.~R\"ulling  for very helpful comments.
This paper is written during a stay at the \'Ecole normale sup\'erieure which is supported by a fellowship within the Post-Doc program of the Deutsche Forschungsgemeinschaft (DFG).
I thank the \'Ecole normale sup\'erieure for its hospitality. 

\section{Zero cycles and the first notch of the coniveau}

\subsection{Notation} 
Let $k$ be a field. By a motive we understand a pair $(X,P)$ with $X$ a smooth projective variety over $k$ and 
$P\in \Hom(X,X)$ a projector in the algebra of correspondences. The correspondences are defined to be  
$$
\Hom(X,Y)=\oplus_i \CH^{\dim X_i}(X_i,Y),
$$
where $X_i$ are the connected components of $X$. Here and in the following we use Chow groups with 
$\Q$ coefficients. Note that we work with effective motives only. 

We simply write $X=(X,id_X)$ for the motive associated with $X$.
The motives form a category $\mathcal{M}_k$ with morphism groups 
$$ 
\Hom((X,P),(Y,Q))= Q \circ \Hom(X,Y) \circ P \subset \Hom(X,Y).
$$
The sum and the product in $\mathcal{M}$ are defined by disjoint union and product:
\begin{align*}
(X,P)\oplus (Y,Q) &= (X\cup Y,P+Q) \\
(X,P)\otimes (Y,Q) &= (X\times Y,P\times Q) 
\end{align*} 

We denote by $\Q(-1)$ the Lefschetz motive, i.e. $\mathbb{P}^1=\Q(0)\oplus \Q(-1)$. 
We set $\Q(a):=\Q(-1)^{\otimes -a}$ for $a<0$ and $\Q(0):=\Spec(k)$.
If $X$ is connected then
$$
\Hom((X,P)\otimes \Q(a),(Y,Q)\otimes \Q(b))= P\circ \CH^{\dim X-a+b}(X\times Y) \circ Q.
$$ 
If $M$ is a motive, we define
$$
\CH^i(M):= \Hom(\Q(-i),M),  \quad  \CH_i(M):=\Hom(M,\Q(-i))
$$
for $i\geq 0$ and $\CH^i(M)=0=\CH_i(M)$ for $i<0$. We have
\begin{equation} \label{Shift}
\CH^i(M\otimes \Q(a))= \CH^{i+a}(M), \quad \CH_i(M\otimes \Q(a))= \CH_{i+a}(M)
\end{equation}
for all $i\geq 0$ and $a\leq 0$. Note that for a motive $M=(X,P)$ with $X$ connected of dimension $n$ 
the equality $\CH_i(M)=\CH^{n-i}(M)$ in general doesn't hold.   

If $k\subset L$ is an extension of fields then $(X,P)\mapsto (X\times_k L,P\times_k L)$ 
defines a functor 
\begin{equation}\label{EL}
\times_k L: \mathcal{M}_k \xr{} \mathcal{M}_L.
\end{equation}

The following Proposition is a consequence of the theory of birational motives \cite{KS} 
due to B.~Kahn and R.~Sujatha. We include the proof for the convenience of the reader.

\begin{proposition} \label{coniveaucondition} 
Let $k$ be a perfect field and $X$ be connected.
\begin{enumerate}
\item[(i)]
A motive $M=(X,P)$ can be written as $M \cong N\otimes \Q(-1)$ with some motive 
$N$ if and only if $\CH_0(M\times_k L)=0$ for some field extension $L$ of the 
function field $k(X)$ of $X$.
\item[(ii)] 
There exists an isomorphism $M \cong N\otimes \Q(a)$ with some motive $N$ 
and $a<0$ if and only if $\CH_i(M\times_k L)=0$ for all $i<-a$ and 
all field extensions $k\subset L$. 
\end{enumerate}
\begin{proof}
(i)\; If $M\cong N\otimes \Q(-1)$ then $M\times_k L\cong (N\times_k L)\otimes \Q(-1)$
and therefore $\CH_0(M\times_k L)=0$ by \ref{Shift}. 

Suppose now that $\CH_0(M\times_k L)=0$. By the same arguments as in \cite[Proposition~1]{BS}
we have 
\begin{equation} \label{PsupportinD}
P \in {\rm image}\left(\CH^{\dim D}(X\times D) \xr{(id \times \imath)_*} \CH^{\dim X}(X\times X)\right)
\end{equation}
for some effective (not necessarily irreducible) Divisor $\imath:D\xr{} X$. For the convenience 
of the reader we recall the proof. It is well-known that 
$$
\CH_0(X\times_k k(X)) \xr{} \CH_0(X\times_k L) 
$$
is injective, and therefore $\CH_0(M\times_k L)=0$ implies $\CH_0(M\times_k k(X))=0$. 
Let $\tau$ be the composite 
$$
\tau: \CH^{\dim X}(X\times X) \xr{} \varinjlim_{U\subset X} \CH^{\dim X}(X\times U)=\CH^{\dim X}(X\times k(X)),
$$
where the limit is over all open subsets $U\subset X$. 
It is easy to see that the equality $0=(P\times_k k(X))\circ \tau(\Delta_X)=\tau(P)$ holds, 
which shows \ref{PsupportinD}.

Let $Y\xr{} D$ be an alteration such that $Y$ is regular (and thus smooth), and denote by 
$f:Y\xr{} D \xr{\imath} X$ the composite. We have $P=(id_X \times f)_*(Z)$ for a suitable 
cycle $Z\in \CH^{\dim Y}(X\times Y)$. Define $Q\in \End(Y)$ by $Q=Z\circ P\circ \Gamma(f)^t$
where $\Gamma(f)^t\in \CH^{\dim X}(Y\times X)$ is the graph of $f$. The equality  
$\Gamma(f)^t\circ Z=P$ implies $Q^2=Q$. It is easy to check that 
$$
(Y,Q)\otimes \Q(-1) \xr{P\circ \Gamma(f)^t} (X,P) \quad (X,P)\xr{Z\circ P} (Y,Q)\otimes \Q(-1)
$$
are inverse to each other, so that $(Y,Q)\otimes \Q(-1)\cong (X,P)$ as claimed. 

(ii)\; By induction on $a$ and using (i).
\end{proof}
\end{proposition}


\subsection{Motives associated with morphism}\label{XZ}  
Let $\pi:X\xr{} Z$ be a finite surjective morphism of degree $d$, where $X$ is connected, smooth and projective, 
but $Z$ may be singular. 
The cycle $X\times_Z X \subset X\times X$ gives a projector $P=1/d \cdot [X\times_Z X]\in \End(X)$ and we
write $(X,\pi):=(X,P)$ for the corresponding motive.

If $\pi:X \xr{} Y$ is a surjective morphism between connected, smooth and 
projective varieties of the same dimension, then the graph $\Gamma(\pi)$ 
of $\pi$ gives morphisms \linebreak 
$\Gamma(\pi)\in \Hom(Y,X)$ and $\Gamma(\pi)^t \in \Hom(X,Y)$. 
Let $d$ be the degree of $\pi$, and $Q$ be a projector; 
since $\Gamma(\pi)^t \circ \Gamma(\pi)=d\cdot id_Y$ the correspondence 
$P=1/d \cdot \Gamma(\pi) \circ Q \circ \Gamma(\pi)^t$ is a projector and $(X,P)\cong (Y,Q)$.
       
\begin{proposition} \label{fromXtoY}
Let $k$ be a perfect field. In the diagram 
$$
\xymatrix
{
X \ar[d]^{\pi} & \\
Z & Y \ar[l]^{f}
}
$$
we assume that $X,Y$ are smooth, connected and projective varieties of the same dimension, the 
morphism $\pi$ is finite and surjective, and $f$ is birational.
The following holds:
\begin{enumerate}
\item[(i)] The motive $(X,\pi)$ is a direct summand in $Y$.
\item[(ii)] If $X=(X,\pi)\oplus N'\otimes \Q(-1)$ for some motive $N'$, then 
$$Y \cong (X,\pi)\oplus N \otimes \Q(-1)$$ for some motive $N$.
\end{enumerate}
\begin{proof}
(i) We write $S$ for the unique irreducible component of $X\times_Z Y$ of dimension $\dim X$. Choose an alteration
$g:W\xr{} S$ with $W$ regular, $W$ is smooth since $k$ is perfect. 

Via $g_1:=pr_1\circ g$ (resp. $g_2:=pr_2\circ g$) the motives $X$, $(X,\pi)$
(resp. $Y$) are direct summands of $W$, we write $P_X,P_{(X,\pi)},P_Y$ for the corresponding projectors. The inclusion
$(X,\pi)$ factors through $Y$ if and only if $P_{(X,\pi)}\circ P_Y = P_Y \circ P_{(X,\pi)}=P_{(X,\pi)}$ in $\End(W)$.
We have
\begin{equation*}
\begin{split}
\deg(g)^2\deg(\pi)^2\cdot P_Y \circ P_{(X,\pi)} &= \Gamma(g_2)  \circ \Gamma(g_2)^t \circ \Gamma(g_1) \circ [X\times_Z X] \circ \Gamma(g_1)^t \\
&= \deg(g)\cdot \Gamma(g_2)  \circ [S] \circ [X\times_Z X] \circ \Gamma(g_1)^t \\
&= \deg(g)\cdot [W\times_Z X] \circ [X\times_Z X] \circ \Gamma(g_1)^t \\
&= \deg(g)\deg(\pi) \cdot \Gamma(g_1) \circ  [X\times_Z X] \circ  [X\times_Z X] \circ \Gamma(g_1)^t \\
&= \deg(g)^2\deg(\pi)^2\cdot P_{(X,\pi)} 
\end{split} 
\end{equation*}
That $P_{(X,\pi)}\circ P_Y=P_{(X,\pi)}$ can be proved in the same way. 
Note that 
$$(X,\pi)\xr{\Gamma(g_1)} W \xr{\Gamma(g_2)^t} Y$$ 
does not depend on the choice of $W$, i.e. $(X,\pi)$ is in a natural way a direct summand in Y. 
Indeed, if $h:W'\xr{} W$ then 
\begin{equation*}
\Gamma(g_2\circ h)^t \circ \Gamma(g_1 \circ h) =  \Gamma(g_2)^t \circ \Gamma(h)^t \circ \Gamma(h) \circ \Gamma(g_1) 
= \Gamma(g_2)^t \circ \Gamma(g_1),     
\end{equation*}
and for another choice $W''$ we may find $W'$ dominating $W$ and $W''$. 

(ii) Write $Y\cong (X,\pi)\oplus M$. Let $L\supset k$ be a field extension, we have 
$
Y\times_k L \cong (X\times_k L,\pi\times_k L) \oplus M\times_k L.
$
The map $S\times_k L \xr{} X\times_k L$ is birational and $X$ is smooth, thus 
$$
\CH_0(S\times_k L) \cong \CH_0(X\times_k L) \cong \CH_0(X\times_k L,\pi\times_k L).
$$
The pushforward $\CH_0(S\times_k L) \xr{} \CH_0(Y\times_k L)$ is surjective, and therefore
$$
\CH_0(Y\times_k L) = \CH_0(X\times_k L,\pi\times_k L)
$$
and $\CH_0(M\times_k L)=0$. 
According to Proposition \ref{coniveaucondition} this shows $M\cong N \otimes \Q(-1)$.
\end{proof}
\end{proposition}

\subsection{Coniveau filtration}
Let $k=\C$, we work with the singular cohomology in rational coefficients 
$H^i(X):=H^i(X,\Q)$ for $i\geq 0$. The coniveau filtration $N^*H^i(X)$ is defined
to be 
$$
N^pH^i(X):= \bigcup_{S} \ker\left(H^i(X)\xr{} H^i(X-S)\right),
$$ 
where $S$ runs through all algebraic subsets (maybe reducible) of codimension $\geq p$. \linebreak
The coniveau filtration is a filtration of Hodge structures and therefore the graduated
pieces $\Gr^p_N:=N^pH^i(X)/N^{p+1}H^i(X)$ inherit a Hodge structure. 

By the work of Arapura and Kang \cite[Theorem~1.1]{AK} the coniveau filtration is 
preserved (up to shift) by pushforwards, exterior products and pullbacks. Using 
resolution of singularities it follows that 
\begin{equation} \label{functorconiveau} 
\Gr^p_N:(X,P) \mapsto {\rm image}(P:\oplus_i \Gr^p_NH^i(X) \xr{} \oplus_i \Gr^p_NH^i(X))
\end{equation} 
is a functor from motives to Hodge structures (for all $p\geq 0$). Note, however, that there
is no Kuenneth formula for $\Gr^p_N$; even for $p=0$ the surjection 
$$
\bigoplus_{s+t=i} \Gr^0_NH^s(X)\otimes \Gr^0_NH^t(Y) \xr{} \Gr^0_NH^i(X\times Y) 
$$
is not injective in general.
For the fiber product with $\PP^1$ we have
$$
N^pH^i(X\times \PP^1)= N^pH^i(X) \oplus N^{p-1}H^{i-2}(X)(-1) 
$$
and therefore
\begin{equation} \label{coniveaushift}
\begin{split}
\Gr^p_N(M\otimes \Q(-1)) &= \Gr^{p-1}_N(M)(-1) \quad \text{if $p>0$} \\
\Gr^0_N(M\otimes \Q(-1)) &= 0 \\
\end{split}
\end{equation} 
for all motives $M$.

\section{Application: the Dwork family and its mirror}

\subsection{} \label{DworkfamilyandGoperation}
Let $k$ be a field. We consider the hypersurfaces $X_\lambda$ in $\PP^n_k$ defined by the equation
\begin{equation} \label{Dworkfamily}
\sum_{i=0}^n X_i^{n+1} + \lambda\cdot X_0\cdots X_n = 0
\end{equation}
with $\lambda \in k$, and we assume that $n+1$ is prime to the characteristic of $k$. 

Let $G\subset (\mu_{n+1})^{n+1}/\Delta(\mu_{n+1})$ ($\Delta(\mu_{n})\cong \mu_{n+1}$ diagonally embedded) 
be the kernel of the character $(\zeta_0,\dots,\zeta_{n+1})\mapsto \zeta_0\cdot \dots \cdot \zeta_{n+1}$, then $G$ acts on
$X_{\lambda}$ in the obvious way. We denote by $\pi:X_{\lambda} \xr{} X_{\lambda}/G $ the quotient map.

\begin{lemma} \label{keylemma}
Let $k$ be a field. We assume that ${\rm char}(k)\nmid n+1$ if ${\rm char}(k)>0$. \linebreak  
If $n\geq 2$ and $X_{\lambda}$ is smooth, then the map  
$$
\CH_0(X_{\lambda}) \xr{} \CH_0(X_{\lambda},\pi)
$$
from section \ref{XZ} is an isomorphism. 
\begin{proof} 
The projector for $(X_{\lambda},\pi) \subset X_{\lambda}$ is $\frac{1}{\abs{G}} \sum_{g\in G} \Gamma(g)$. Therefore 
the statement is equivalent to 
$$
\sum_{g\in G} g_*(a) = \abs{G} \cdot a
$$
for every $a\in \CH_0(X_{\lambda})$.

1.~case: $k=\C$. For $n=2$  the quotient map $\pi:X_{\lambda} \xr{} X_{\lambda}/G$ is an isogeny of elliptic curves, and therefore the statement is true.   
 
Consider $\mu_{n+1}\cong H\subset G$ with $\zeta \mapsto (\zeta,\zeta^{-1},1,\dots,1)$, and $\tau\in \Aut(X_{\lambda})$ defined by 
$\tau^*(X_0)=X_1, \tau^*(X_1)=X_0,$ and $\tau^*(X_i)=X_i$ otherwise. We have $H\rtimes \Z/2\cdot \tau \subset \Aut(X_{\lambda})$ and 
claim that $X_{\lambda}/(H\rtimes \Z/2\cdot \tau)$ is rational. Indeed, for the open set $U_{\lambda}=\{ X_n\neq 0 \}\subset X_{\lambda}$
we compute 
\begin{equation*}
\begin{split}
U_{\lambda}/(H\rtimes \Z/2\tau)\cong \Spec(k[\sigma_1,x_2,\dots,x_{n-1},v]/I)\cong \Spec(k[x_2,\dots,x_{n-1},v]),\\
\end{split}
\end{equation*}
with $I=(\sigma_1+x_2^{n+1}+\dots+x_{n-1}^{n+1}+\lambda\cdot v \cdot x_2\dots x_{n-1})$. Here, the coordinates
are defined to be $x_i:=X_i/X_{n}$, $v=x_0\cdot x_1$, and $\sigma_1= x_0^{n+1}+x_1^{n+1}$. Since rational varieties are 
rationally chain connected we conclude that
\begin{equation}\label{identityHxtau}
\sum_{g\in H\rtimes \Z/2\tau} \Gamma(g)\cdot a = 2(n+1) \cdot  \deg(a) \cdot [p] 
\end{equation}
for every $a\in \CH_{0}(X_{\lambda})$ and some closed point $p\in X_{\lambda} \cap \{X_0=X_1=0\}$ ($p$ exists since $n>2$).
 
Next, if $\zeta\in \mu_{n+1}\cong H$ then $(\zeta,\tau)\in H\rtimes \Z/2\cdot \tau$ has order $2$, and we consider the
quotient $q:X_{\lambda}\xr{}X_{\lambda}/(\zeta,\tau)$ by the action of $(\zeta,\tau)$. We claim that $X_{\lambda}/(\zeta,\tau)$
is rationally chain connected.

The fixpoint set $F$ is 
\begin{equation*}
\begin{split}
&F=\{ X_0-\zeta X_1 = 0\} \quad \text{if $n$ is odd}, \\
&F=\{ X_0-\zeta X_1 = 0\} \cup \{ [1:-\zeta^{-1}:0:\dots :0]\} \quad \text{if $n$ is even}. 
\end{split}
\end{equation*}
Let $H=\{ X_0-\zeta X_1 = 0\}\subset F$ be the hyperplane section. One verifies that $H$ is smooth
if and only if $X_{\lambda}$ is smooth, and for every point $x\in H$ there are coordinates $y,x_1,\dots,x_{n-2}$
such that $y$ is a local equation for $H$ with $(\zeta,\tau)^* y =-y$ and the $x_i$ are invariant. 
Thus $y^2,x_1,\dots,x_{n-2}$ are local coordinates for the quotient which is therefore smooth in the points $q(H)$.
So that $X_{\lambda}/(\zeta,\tau)$ is smooth if $n$ is odd, and $X_{\lambda}/(\zeta,\tau)$ has an isolated quotient 
singularity in $q([1:-\zeta^{-1}:0:\dots :0])$ if $n$ is even.

In both cases, $2K_{X_{\lambda}/(\zeta,\tau)}$ is Cartier and $2K_{X_{\lambda}/(\zeta,\tau)}\cong \OO(-q(H))$ (the isomorphism
comes from an invariant form in $H^0(X_{\lambda},\omega_{X_{\lambda}}^{\otimes 2})=H^0(X_{\lambda},\omega_{X_{\lambda}}^{\otimes 2})^{(\zeta,\tau)}$). 
We have $q^*(\OO(q(H)))=\OO(2H)$ and therefore $\OO(q(H))$ is ample. If $n$ is odd then the Theorem of Campana, Koll\'ar, Miyaoka, Mori (\cite{C},\cite{KMM}) implies 
that $X_{\lambda}/(\zeta,\tau)$ is rationally chain connected. If $n$ is even, then $X_{\lambda}/(\zeta,\tau)$ is a $\Q$-Fano 
variety with log terminal singularities and we may use the Theorem of Zhang \cite{Z} to  prove the claim.

We conclude that 
\begin{equation} \label{identityzetaxtau}
a + \Gamma((\zeta,\tau))(a) = 2\deg(a)[p]
\end{equation}
for every $a\in \CH_0(X_{\lambda})$ and $p\in X_{\lambda} \cap \{X_0=X_1=0\}$. 
Using \ref{identityHxtau} and \ref{identityzetaxtau} we see
\begin{equation}
\begin{split} \label{identity}
\sum_{g\in H} \Gamma(g)(a) &= \sum_{g\in H\rtimes \Z/2\tau} \Gamma(g)(a) - \sum_{\zeta\in \mu_{n+1}} \Gamma((\zeta,\tau))(a) \\
                        &= 2(n+1)\deg(a)[p] - \sum_{\zeta\in \mu_{n+1}} (2\deg(a)[p]-a)= (n+1) a.
\end{split}
\end{equation}
Of course, for the subgroups $\mu_{n+1}\cong H_i \subset G$ defined by 
$\zeta\mapsto (1,\dots,1,\zeta,\zeta^{-1},1,\dots,1)$ where $\zeta$ is put in the $i$-th position, the same conclusion
\ref{identity} holds. Now, the equality 
$$
\sum_{g\in G} \Gamma(g) = \left( \sum_{g\in H_0} \Gamma(g) \right) \circ \dots \circ \left( \sum_{g\in H_{n-2}} \Gamma(g) \right)
$$  
proves the claim.

2.~case: ${\rm char}(k)=0$. It is a well-known fact that if $k_0\subset k$ is a subfield and 
$X=X_0\times_{k_0}k$ then the pullback map
\begin{equation} \label{CH0k0k}
\CH_0(X_0) \xr{} \CH_0(X)
\end{equation}
is injective (without the assumption on ${\rm char}(k)$). 
The variety $X_{\lambda}$ is defined over 
$\Q(\lambda)\subset k$, and every zero cycle can be defined over a subfield $k_0\subset k$
which is finitely generated over $\Q(\lambda)$. By fixing an embedding $\sigma:k_0\xr{} \C$, 
we reduce to the case $k=\C$. 

3.~case: ${\rm char}(k)=p\neq 0$. Again, since \ref{CH0k0k} is injective, 
we may assume that $k$ is algebraically closed. Let $W$ 
be the Witt vectors of $k$; $W$ is a complete discrete valuation ring with residue field 
$k$ and quotient field $K$ with ${\rm char}(K)=0$. Choose a lift $\tilde{\lambda}\in W$ 
of $\lambda$, and let $X_{\lambda,W}\subset \PP^n_W$ be the variety 
$\sum_{i=0}^n X_i^{n+1}+\tilde{\lambda} X_0\cdots X_n=0$. The specialization map 
$$
sp: \CH_0(X_{\lambda,W}\otimes_W K) \xr{} \CH_0(X_{\lambda})
$$   
from \cite[\textsection20.3]{Fulton} is surjective, because $W$ is complete (and therefore 
$X_{\lambda,W}(W)\xr{} X_{\lambda}(k)$ is surjective). 
Since ${\rm char}(k)\nmid n+1$ we have 
$$
\mu_{n+1}(k) \xl{\cong} \mu_{n+1}(W) \xr{\cong} \mu_{n+1}(K),
$$
and the same statement holds for $G$.
Now the compatibility 
of $sp$ with pushforwards \cite[Proposition~20.3]{Fulton} proves the claim. 
\end{proof}
\end{lemma}

\begin{thm}\label{thm}
Let $k$ be a perfect field. We assume ${\rm char}(k)\nmid n+1$ if ${\rm char}(k)>0$. 
Let $X_{\lambda}$ be a smooth member of the Dwork family for $n\geq 2$.   
If $\pi:X_{\lambda} \xr{} X_{\lambda}/G$ 
is the quotient of the $G$-action (see \ref{DworkfamilyandGoperation}) and
$Y_{\lambda} \xr {} X_{\lambda}/G$ is a resolution of singularities, then 
$$
X_{\lambda}\cong (X_{\lambda},\pi) \oplus N'_{\lambda}\otimes \Q(-1), \quad
Y_{\lambda}\cong (X_{\lambda},\pi) \oplus N_{\lambda}\otimes \Q(-1)
$$
for some motives $N'_{\lambda}$ and $N_{\lambda}$.
\begin{proof}
By construction of $(X_{\lambda},\pi)$ in \ref{XZ} we have 
$X_{\lambda}\cong (X_{\lambda},\pi)\oplus M_{\lambda}$ with some motive $M_{\lambda}$.
In view of Lemma \ref{keylemma} we know that 
$$
\CH_0(X_{\lambda}\times_k L)= \CH_0((X_{\lambda}\times_k L,\pi\times_k L))= 
\CH_0((X_{\lambda},\pi)\times_k L)
$$
for all field extensions $k\subset L$, and thus $\CH_0(M_{\lambda}\times_k L)=0$. 
Proposition \ref{coniveaucondition} implies that $M_{\lambda}\cong N'_{\lambda}\otimes \Q(-1)$
for some $N'_{\lambda}$, and Proposition \ref{fromXtoY} proves the claim. 
\end{proof}
\end{thm}

\begin{korollar} \label{consequences} Under the assumptions of Theorem \ref{thm}. 
\begin{enumerate}
\item[(i)] If $k=\C$ then there is an isomorphism of Hodge structures 
$$
\Gr^0_NH^*(X_{\lambda},\Q)\cong \Gr^0_NH^*(Y_{\lambda},\Q). 
$$
\item[(ii)] If $k=\mathbb{F}_q$, the finite field with $q$ elements, then for all $m\geq 1:$
$$
\#X_{\lambda}(\mathbb{F}_{q^m}) = \#Y_{\lambda}(\mathbb{F}_{q^m}) \quad \text{modulo $q^{m}$.} 
$$
\end{enumerate}
\begin{proof}
(i)\; By \ref{coniveaushift}. 

(ii) If $N$ is a motive (over $\mathbb{F}_q$) then the eigenvalues of the Frobenius acting on 
$H^*_{\text{\'et}}(N\otimes \Q(-1))=H^*_{\text{\'et}}(N)\otimes \Q_l(-1)$ lie in 
$q\cdot \bar{\Z}$. Now the claim follows from Grothendieck's trace formula. 
\end{proof}
\end{korollar}

\section{Conjectures}

\subsection{}
For $k=\C$ the Hodge structure of a
variety $X$ can be recovered from the associated motive. For an effective
motive $N$ we know $h^{i,0}(N\otimes \Q(-1))=0$ for all $i$, so that 
if
\begin{equation}\label{statementmotives}
X\cong X' \oplus N'\otimes \Q(-1), \quad
Y\cong X' \oplus N\otimes \Q(-1),
\end{equation}
then $h^{i,0}(X)=h^{i,0}(Y)$.

\subsection{} 
Now consider the setting 
$$
\xymatrix
{
X \ar[d]^{\pi} & \\
X/G & Y \ar[l]
}
$$
where $X$ is a smooth projective variety with an action of a finite group $G$, 
and $Y$ is a resolution of singularities of the quotient $X/G$. Since $X/G$
has only quotient singularities we know that $H^i(Y,\OO_Y)=H^i(X/G,\OO_{X/G})$ for
all $i$. The map $\OO_{X/G}\xr{} \pi_*\OO_X$ is split by $\frac{1}{\abs{G}}\sum_{g\in G}g^*$ and
we obtain 
$$
H^i(X,\OO_{X})^G=H^i(X/G,\OO_{X/G})=H^i(Y,\OO_Y), \quad \text{for all $i$.}
$$  
Therefore \ref{statementmotives} can only be expected if the following holds:
\begin{equation}\label{conditioncohomology}
H^i(X,\OO_{X})=H^i(X,\OO_{X})^G. 
\end{equation}

\subsection{} On the other hand, the Bloch conjectures on a filtration of the Chow group 
of zero cycles which is controlled by the Hodge structure (see \cite[\S 23.2]{V} for a
precise statement) predict
$$
\pi:\CH_0(X)\xr{\cong} \CH_0(X/G)=\CH_0(X,\pi)
$$
whenever \ref{conditioncohomology} holds, and thus $X\cong (X,\pi)\oplus N'\otimes \Q(-1)$
(in the notation of \ref{XZ}). Now, Proposition \ref{fromXtoY} yields \ref{statementmotives}
with $X'=(X,\pi)$. So that the Bloch conjectures imply the following conjecture.

\begin{conjecture}\label{conj}
Let $X$ be a smooth projective variety over a field $k$ of ${\rm char}(k)=0$, and
let $G$ be a finite group acting on $X$ with $H^i(X,\OO_{X})=H^i(X,\OO_{X})^G$ for all 
$i$. If $Y$ is a resolution of singularities of the quotient $\pi:X\xr{} X/G$, then there
are (effective) motives $N,N'$ such that    
$$
X\cong (X,\pi) \oplus N'\otimes \Q(-1), \quad
Y\cong (X,\pi) \oplus N\otimes \Q(-1).
$$
\end{conjecture}

Unfortunately little is known concerning the Bloch conjecture.

\subsection{} 
Let us consider monomial deformations of the degree $d$ Fermat hypersurface in $\PP^n$, and
$G\subset \mu_d^{n+1}$.  The condition \ref{conditioncohomology} holds only for 
$d\leq n+1$ (or $G=\{1\}$). Therefore there is no generalisation of Theorem \ref{thm}
to degree $d>n+1$. In the case $d<n+1$, $X_{\lambda}$ and $Y_{\lambda}$ are rationally connected,
and thus $\CH_0(X_{\lambda})=\Q=\CH_0(Y_{\lambda})$. We obtain 
$$
X_{\lambda} \cong \Q \oplus N'\otimes \Q(-1), \quad
Y_{\lambda} \cong \Q \oplus N\otimes \Q(-1).
$$
from Proposition \ref{coniveaucondition}. 

\subsection{}
For a finite field $k=\mathbb{F}_q$ we don't know the correct assumptions 
for Conjecture \ref{conj}. However, the assertion implies a 
congruence formula for the number of rational points:
\begin{equation}\label{congXY}
\#X(\mathbb{F}_q) \equiv \#Y(\mathbb{F}_q) \mod q.
\end{equation}    

The work of Fu and Wan provides a congruence formula for the number 
of rational points of $X$ and $X/G$. 

\begin{thm}[\cite{FW}]\label{thmFW}
Let $X$ be a smooth projective
variety over the finite field $\mathbb{F}_q$. Suppose $X$ has a smooth projective lifting
$\mathcal{X}$ over the Witt ring $W=W(\mathbb{F}_q)$ such that the $W$-modules
$H^r(\mathcal{X},\Omega_{\mathcal{X}/W}^s)$ are free. Let $G$ be a finite group of
$W$-automorphisms acting on $X$. Suppose $G$ acts
trivially on $H^i(\mathcal{X}, {\OO}_{\mathcal{X}})$ for all $i$. Then for any
natural number $k$, we have the congruence
$$\# X(\mathbb{F}_{q^k})\equiv \# (X/G) (\mathbb{F}_{q^k})~({\rm mod}~q^k).$$
\end{thm}

By extending the theory of Witt vector cohomology to singular varieties, Berthelot, Bloch 
and Esnault were able to prove the following theorem.

\begin{thm}[{\cite[Corollary~6.12]{BBE}}]\label{thmBBE}
Let $X$ be a proper scheme over $\mathbb{F}_q$, and $G$ a finite group acting
on $X$ so that each orbit is contained in an affine open subset of
$X$. If $\abs{G}$ is prime to $p$, and if the action of $G$ on $H^i(X,
\OO_X)$ is trivial for all $i$, then
$$
\#X(\mathbb{F}_q) \equiv \#(X/G)(\mathbb{F}_q) \mod q. 
$$ 
\end{thm}

In the next section we prove that 
$$
\#Y(\mathbb{F}_q) \equiv \#(X/G)(\mathbb{F}_q) \mod q
$$
if $X$ is a smooth projective variety and $Y$ is a resolution of 
singularities of $X/G$. So that with the assumptions of \ref{thmFW} or \ref{thmBBE}  
we obtain the congruence formula \ref{congXY}.  

\section{Congruence formula}

\subsection{}
In this section we fix a finite field $k=\mathbb{F}_q$ of characteristic $p$, and an 
algebraic closure $\bar{k}$ of $k$. For a separated scheme $X$ of finite type over $k$ 
we work with the \'etale cohomology groups $H^i(X_{\bar{k}},\Ql)$ (resp. \'etale cohomology 
groups with support $H^i_{Z_{\bar{k}}}(X_{\bar{k}},\Ql)$ for $Z\subset X$) with $\ell$ a 
prime number $\not=p$, and $X_{\bar{k}}=X\times_k \bar{k}$. They are finite dimensional 
$\Ql$-vector spaces, and the Galois group $G_k={\rm Gal}(\bar{k}/k)$ acts continuously 
on them. 

We denote by $F\in G_k$ the geometric Frobenius, $F$ acts on 
$H^i(X_{\bar{k}},\Ql)$ with eigenvalues that are algebraic integers. If $X$ is proper
then we have Grothendieck's trace formula 
$$
\#X(\mathbb{F}_q)=\sum_i (-1)^i {\rm Tr}(F,H^i(X_{\bar{k}},\Ql)).
$$

\subsection{}
Let $V$ be a finite dimensional $\Ql$-vector space, and $F:V\xr{} V$ a linear map.
Fix an algebraic closure $\bar{\Q}_{\ell}$ of $\Ql$. 
The vector space $\bar{V}=V\otimes_{\Ql} \bar{\Q}_{\ell}$ decomposes into the generalised 
eigenspaces of $F$: 
$$
\bar{V}=\bigoplus_{\lambda \in \bar{\Q}_{\ell}} \bar{V}_{\lambda},
$$
i.e.~$\bar{V}_{\lambda}$ is the maximal subspace such that $F$ acts with eigenvalue $\lambda$.
For every $g\in G_{\Ql}={\rm Gal}(\bar{\Q}_{\ell}/\Ql)$ we get $g(\bar{V}_{\lambda})=\bar{V}_{g(\lambda)}$, and we obtain a decomposition 
$$
V= \bigoplus_{\lambda\in G_{\Q_{\ell}}\backslash \bar{\Q}_{\ell}} \left( \bigoplus_{\lambda'\in G_{\Ql}\cdot \lambda} \bar{V}_{\lambda}\right)^{G_{\Ql}},
$$
where $\lambda$ runs through all orbits of $G_{\Q_{\ell}}$ in $\bar{\Q}_{\ell}$, and $\lambda'$
through all conjugates of $\lambda$. We write 
$$
V_{\lambda}=\left( \bigoplus_{\lambda'\in G_{\Ql}\cdot \lambda} \bar{V}_{\lambda}\right)^{G_{\Ql}}, \quad V=\bigoplus_{\lambda\in G_{\Q_{\ell}}\backslash \bar{\Q}_{\ell}} V_{\lambda}.
$$ 
Let $W$ be another finite dimensional $\Ql$-vector space with a linear operation $F:W\xr{} W$.
If $\phi:V\xr{} W$ is a linear map which commutes with the action of $F$ then 
$$
\phi(V_{\lambda})\subset W_{\lambda}
$$
for every $\lambda\in G_{\Q_{\ell}}\backslash \bar{\Q}_{\ell}$.

Now, we fix an integer $q\in \Z$, and we assume that all eigenvalues of $F$ are algebraic integers, 
i.e.~$V_{\lambda}=0$ if $\lambda\not\in G_{\Q_{\ell}}\backslash \bar{\Z}$ where $\bar{\Z}\subset\bar{\Q}_{\ell}$ is the integral closure of $\Z$ in $\bar{\Q}_{\ell}$.
We note that the subset $q\bar{\Z}\subset \bar{\Z}\subset \bar{\Q}_{\ell}$ has an induced action
by $G_{\Q_{\ell}}$, and we define the slope $<1$ resp. slope $\geq 1$ part of $V$ to be
$$
V^{<1}:=\bigoplus_{\lambda\not\in G_{\Q_{\ell}}\backslash q\bar{\Z}} V_{\lambda}, \quad V^{\geq 1}:=\bigoplus_{\lambda\in G_{\Q_{\ell}}\backslash q\bar{\Z}} V_{\lambda}. 
$$ 
We obtain a decomposition 
$$
V=V^{<1}\oplus V^{\geq 1}
$$
with $F$ action on $V^{<1}$ and $V^{\geq 1}$, and the decomposition is functorial for linear maps
that commute with the $F$-operation.

For \'etale cohomology and $q=\abs{k}$ we thus get for all $i$ a functorial decomposition 
\begin{equation*}
\begin{split}
H^i(X_{\bar{k}},\Ql)&=H^i(X_{\bar{k}},\Ql)^{<1}\oplus H^i(X_{\bar{k}},\Ql)^{\geq 1}, \\
\text{resp.} \quad H^i_{Z_{\bar{k}}}(X_{\bar{k}},\Ql)&= H^i_{Z_{\bar{k}}}(X_{\bar{k}},\Ql)^{<1}\oplus 
H^i_{Z_{\bar{k}}}(X_{\bar{k}},\Ql)^{\geq 1}. 
\end{split}
\end{equation*}

\begin{lemma}\label{lemmaHelene}
If $X$ is smooth and $Z\subset X$ is a closed subset of codimension $\geq 1$ then 
$$H^i_{Z_{\bar{k}}}(X_{\bar{k}},\Ql)^{<1}=0 \quad \text{for all $i$.}$$
\begin{proof}
\cite[Lemma~2.1]{E1},\cite[\S 2.1]{E2}.
\end{proof}
\end{lemma}

In other words all eigenvalues of the Frobenius on $H^i_{Z_{\bar{k}}}(X_{\bar{k}},\Ql)$ lie in $q\cdot \bar{\Z}$.
It is not difficult to extend this lemma to the case when $X$ has quotient singularities. 

\begin{lemma}\label{lemmaquotient}
Let $X$ be a smooth and quasi-projective variety, and let $G$ be a finite group acting on $X$. 
If $\pi:X\xr{} X/G$ is the quotient and $Z\subset X/G$ is a closed subset of codimension $\geq 1$ 
then $$H^i_{Z_{\bar{k}}}((X/G)_{\bar{k}},\Ql)^{<1}=0 \quad \text{for all $i$.}$$
\begin{proof}
We write $Y=\pi^{-1}(Z)$ which is a closed subset of $X$ of codimension $\geq 1$. 
Note that $(X/G)_{\bar{k}}=X_{\bar{k}}/G$, i.e.~$\pi_{\bar{k}}:X_{\bar{k}}\xr{} (X/G)_{\bar{k}}$ is the quotient for the $G$ action on $X_{\bar{k}}$. 

The composite of $\Q_{\ell}\subset \pi_{\bar{k}*}\Q_{\ell}$ with  
$$
\sum_{g\in G}g^*:\pi_{\bar{k}*}\Q_{\ell}\xr{} \Q_{\ell}
$$
is multiplication by $\abs{G}$. Since 
$$
H^i_{Z_{\bar{k}}}((X/G)_{\bar{k}},\pi_{\bar{k}*}\Q_{\ell})=H^i_{Y}(X_{\bar{k}},\Ql)
$$
we get 
$$
H^i_{Z_{\bar{k}}}((X/G)_{\bar{k}},\Q_{\ell})\cong H^i_{Y}(X_{\bar{k}},\Ql)^G,
$$
and this map is compatible with the Frobenius action. Now, Lemma \ref{lemmaHelene} implies
the statement.
\end{proof}
\end{lemma}

\begin{thm}\label{thm2}
Let $X$ be a smooth projective $\mathbb{F}_q$-variety with an action of a finite group $G$.
Let $\pi:X\xr{} X/G$ be the quotient, and $f:Y\xr{} X/G$ be a birational map, where $Y$ is 
a smooth projective variety. Then
$$
\#Y(\mathbb{F}_q)=\#X/G(\mathbb{F}_q) \quad \mod q.
$$  
\begin{proof}
Let $U$ be an open (dense) subset of $X/G$ such that $f^{-1}(U)\xr{\cong} U$ is an isomorphism.
Write $Z=(X/G)\backslash U$ and $Z'=Y\backslash f^{-1}(U)$. We consider the map of long exact 
sequences
$$
\xymatrix{
 \ar[r] & H^i_{Z'_{\bar{k}}}(Y_{\bar{k}},\Ql) \ar[r] & H^i(Y_{\bar{k}},\Ql) \ar[r] & H^i(f^{-1}(U)_{\bar{k}},\Ql) \ar[r] & 
\\
\ar[r] & H^i_{Z_{\bar{k}}}((X/G)_{\bar{k}},\Ql) \ar[r]\ar[u] & H^i((X/G)_{\bar{k}},\Ql) \ar[r]\ar[u] & H^i(U_{\bar{k}},\Ql) \ar[u] \ar[r] &
}
$$
Here all maps commute with the action of the Frobenius. By using Lemma \ref{lemmaquotient} we get 
$$
\xymatrix
{
H^i(Y_{\bar{k}},\Ql)^{<1} \ar[r]^{\cong} & H^i(f^{-1}(U)_{\bar{k}},\Ql)^{<1} \\
H^i((X/G)_{\bar{k}},\Ql)^{<1} \ar[r]^{\cong}\ar[u] & H^i(U_{\bar{k}},\Ql)^{<1}. \ar[u]^{=}  &
}
$$
This implies 
$$
H^i(Y_{\bar{k}},\Ql)^{<1} \cong H^i((X/G)_{\bar{k}},\Ql)^{<1} \quad \text{for all $i$.}
$$
With Grothendieck's trace formula we obtain
$$
\#Y(\mathbb{F}_q)-\#X/G(\mathbb{F}_q)= \sum_i (-1)^i 
\left({\rm Tr}(F,H^i(Y_{\bar{k}},\Ql))^{\geq 1}- {\rm Tr}(F,H^i((X/G)_{\bar{k}},\Ql))^{\geq 1}\right).
$$
The right-hand side is a number in $\Z\cap q\bar{\Z}=q\Z$, which proves the congruence. 
\end{proof}
\end{thm}

\begin{remark}
It seems that the fibre $f^{-1}(x)$ of a point $x\in (X/G)(\mathbb{F}_q)$ satisfies the 
congruence 
$$
\#f^{-1}(x)(\mathbb{F}_q) =1 \quad \mod q. 
$$
Of course this would imply the statement of Theorem \ref{thm2}.
\end{remark}

\end{document}